\documentclass[11pt
]{article}
\usepackage{graphicx,amssymb,amsmath}
 \usepackage{tikz,pgf}
 \usetikzlibrary{matrix, arrows, decorations.pathmorphing}

 \usepackage[usenames,dvipsnames]{pstricks}
\usepackage{epsfig
,amscd
}

\numberwithin{equation}{section}
\allowdisplaybreaks
\begin{document}
\title{\large\bf       Detecting Potential Instabilities  of Numerical  Algorithms  }
\author{Yao Yang\\
{\small
  }\\  \\
{\small }\\
{\small }\\[3mm]}

\date{}
\maketitle

\noindent \textbf{\small Abstract }It has been the standard teaching of today that backward stability analysis  is taught as absolute, just as in Newtonian physics time is taught as absolute time. We will prove that it is not true in general.  It depends on algorithms. We will prove  that  forward and mixed stability analysis are absolutely invalid stability analysis in the sense that they have absolutely wrong reference points for detecting huge element growth of any algorithms (if any),  even an "ideal" or "desirable"  backward stability analysis is not so "ideal" or "desirable" in general.  Any of forward stable, backward stable and mixed stable algorithms as in Demmel, Kahan and Parlett and others' papers and text books  \cite{DM}, and \cite{Higham}  may not   be really stable at all because they may fail to detect and expose any potential instabilities of the algorithm in corresponding  stability analysis.   Therefore, it is impossible to prove an algorithm is stable according to the standard teaching of today, just as it is impossible to prove a  mathematical equation(s)(Maxwell's,) is a law of physics according to the standard teaching in Newtonian physics.

 \vskip 1em


 \vskip 1em
\section{Introduction} Given an interesting   numerical algorithm, mathematicians   would like to prove that the  algorithm is stable.   Berkeley's  math professors Kahan. Demmel and Parlett  have been using forward, backward and mixed stability analysis for proving their algorithms are stable in three different stability analysis as in \cite{DK}, \cite{DV} and \cite{FP}.    The standard teaching in today's text books  teaches  backward stabilitt analysis is "desirable" or "ideal"  as in \cite{DM},p.5 or \cite{DP},p.872 and it is taught as absolute stability analysis just as in classical physics in which time is taught as absolute time.  Math professors such as Demmel,  Kahan and Parlett at Berkeley teach  their students if one can derive a backward error bound for an algorithm,  then she/he can claim that  algorithm is proven backward stable and it is "ideal" case scenario.  In all of their accurate singular value computing papers such as award winning paper \cite{DK} and \cite{FP}, they proved their algorithms are forward stable or mixed stable.  Are these forward stable,  backward stable and  mixed stable  algorithms really stable ?  We will prove that the answer is negative in general. We will prove backward stability analysis is not absolute. It depends on algorithms, see Theorem 3 and discussion right before it. Therefore, it is impossible to prove that an algorithm is stable according to standard teaching of today just as it is impossible to prove mathematical equations(e.g. Maxwell) is a law of physics according to standard teaching in classical physics.

\vskip 1em
\noindent  In 1993, the author as a student at Berkeley observed  in  \cite{YY} that given an extremely simple and straightforward algorithm dqd or dqds, his math professors  Kahan, Demmel and Parlett at Berkeley could  not determine if the algorithm has potential instabilities or not in \cite{FP},p.2 which predicts  that dqds has no  potential instabilities at all with any shifts. Just take a look at dqds one can see easily that with a right choice of shift, cancellation can occur and huge element growth can occur.   `   Any testing with dqds would indicate that dqds has infinitely potential instabilities due to excessive element growth with general input vectors.  It follows that in 1994  they could not be sure if dqds has potential instabilities or not in any numerical computing even after they had derived  tiny relative error bound.  Any perturbation result can not help them a bit with the matter because it obviously can not determine if an algorithm has potential instabilities or not either.  It means that in 1994  they could not prove that dqds is really stable based on their mixed stability analysis despite they have proved it is mixed stable(as opposed to forward stable or backward stable).   Whether or not an algorithm has potential instabilities or not can not be determined by deriving relative mixed  error bound with or without any perturbation result. This statement could be generalized to any error analysis as follows: Since one of three stability  analyses, their mixed stability analysis fails totally to be able to detect potential instabilities of an algorithm, dqds viewed as general vector iterations. It is possible that the other two stability analyses, backward stability analysis and forward stability analysis may  also fail to detect potential instabilities due to excessive element growth unless proven otherwise. In 1994, the author asked them if they had such a proof, the answer was negative(this paper will prove such a proof does not exist). Otherwise, in 1994  they would have had a proof to show the author that any of the other two stability analyses is  able to detect potential instabilities.  Since 1994, they probably have been trying to give  such a proof. This paper will prove that such a proof does not exist, see Theorem 2, 3.   Therefore, the student proved that they were not sure if any algorithm has potential instabilities or not in 1994. It implies that his three famous Berkeley math professors were not sure they could compute any digit correctly in any numerical computing( a testing of the author in \cite{YY} shows that no correct digit can be expected if instabilities occur due to excessive element growth). This paper will prove that even "ideal" or "desirable" backward stability analysis (see  \cite{DP},p.872, or \cite{DM},p.5) is not able to detect excessive element growth in general.
\vskip 1em
\noindent There was a second  way to prove in the author's  Ph.D thesis \cite{YY}  that the three famous Berkeley math professors Demmel, Kahan and Parlett were not sure in 1994 they could compute any digit correctly in any numerical computing in theory.   Given an extremely simple and straightforward algorithm dqd or dqds, the author observed and conjectured at Berkeley in his thesis \cite{YY},p.24 that they would never be able to detect potential instabilities of dqds.(I am sure it is true, easy to check) in their standard perspectives.  As a consequence,  their paper \cite{FP} contradicts to a student's Ph.D thesis in which his backward stability analysis  predicts dqds as a general vector iteration has infinitely many potential instabilities due to excessive element growth, being  consistent with their own testing and the author's and consistent with their own teaching but contradicting completely to their mixed stability analysis in \cite{FP},p.2 which predicts that dqds has no potential instabilities at all, what so ever. That is, the author proved that their paper contradicts to their own teaching and to testing with dqds viewed as general vector iterations. Again  they did not have a proof that any of the other two stability analyses, backward stability analysis and forward stability analysis is able to detect potential instabilities. It follows that they were not sure in 1994 any algorithm has potential instabilities or not. This contradicts to their claims that singular values can be computed highly accurately, including tiny ones in \cite{FP} and award winning paper \cite{DK}. Therefore,  the author proved in 1994  automatically right  after he solved the open problem of accurate singular value computing for deriving a backward stability analysis for dqds that they were not sure that  they could compute any digit correctly in any numerical computing in theory.(they did not have a proof that  any of  the other two error analyses, forward and backward error analysis was able to detect potential instabilities )

\vskip 1em
 \noindent  The purpose of listing contradictions above is to prove they have not understood accurate singular values yet up to today.  To the best knowledge of the author, these contradictions still  exist today. The author's  conjecture in \cite{YY} is still true today that they would never be able to detect and expose potential instabilities of dqd or dqds. And the author has shown in his Ph.D thesis \cite{YY} that a backward error analysis can expose all  infinitely many potential instabilities of dqd and dqds but they have neve been able to derive one in terms of their perspectives.   Therefore, they still can not determine if dqds has potential instabilities or not in any numerical computing. And they still do not have a proof that any of the other two stability analyses is able to detect potential instabilities. This paper will prove such a proof does not exist in general.  In mathematics, the self contradictions of theirs usually proves that the three Berkeley math professors have not understood the subject yet.

\vskip 1em

 \noindent  This paper presents a new perspective that is different from  the two proofs listed above as in his Ph.D thesis \cite{YY}  that in all of their accurate singular value computing papers \cite{DK} and \cite{FP} and numerous other scientific computing papers of theirs, the famous  Berkeley math professors Kahan, Demmel and Parlett  have been using forward stability analysis and mixed stability analysis, we will prove forward stable algorithm and mixed stable algorithm may still have potential instabilities due to excessive element growth in general just  as dqds in \cite{FP},p.2  viewed as general vector iterations. More surprisingly, we will prove even an algorithm  that has been  proven backward stable based  on "ideal" or "desirable"  backward stability analysis  may still  have potential instabilities due to excessive element growth. That is, no error analysis is able to detect potential instabilities due to excessive element growth. In particular,  backward stability analysis is not absolute. It depends on algorithms, contrary to the standard teaching of today.  It follows that it is impossible to prove that an algorithm is stable according to the standard teaching of today as in their text book  \cite{DM},p.5 at Berkeley. According to Wilkinson, deriving an error bound itself is usually the least important. Much more important is to expose(detect) potential instabilities, if any, of an algorithm such as dqd or dqds. The following Theorem 1 is consistent with Wilkinson's perspectives as in \cite{GG},p.64.

 \vskip 1em

\noindent\textbf {Theorem 1} Any true stability analysis must satisfy a  necessary condition that it is able to  detect potential instabilities of an algorithm  due to excessive element growth. Excessive element growth can cause potential instabilities due to potential loss of information in any numerical computing and is too important to be ignored in a true stability analysis.

\vskip 1em

\noindent Proof: Otherwise the stability proof or analysis of an algorithm based on the stability analysis will be  inconsistent with computing experience. Excessive element growth can cause potential instabilities of an algorithm due to potential loss of information and is too important to be ignored.  Proof done.

\vskip 1em
\noindent   An extremely simple and straightforward dqd or  dqds is an example that  they have never been able to prove it is stable in any numerical computing, in particular, in singular value computing based on  their mixed stability analysis for a simple reason: they can not determine if dqds has potential instabilities or not based on their mixed stability analysis or any mathematical analysis of theirs. Their mixed stability analysis for dqds viewed as general vector iterations such as \cite{DM},p;.247 or originally in \cite{FP},p.2  have  proved  that dqds has no potential instabilities what so ever. Any simple testing  can verify that dqds viewed as general vector iterations has infinitely many potential instabilities due to excessive element growth.
\vskip 1em
\noindent One can find forward stability analysis and mixed stability analysis in numerous papers of Demmel, Kahan and Parlett such as in \cite{DK}, \cite{FP}, etc and in  all of their singular value computing papers. It follows that  Kahan and his team Demmel and Parlett and others are not sure if any singular value computing has potential instabilities or not.  See Higham \cite{Higham}.

\vskip 1em
 \noindent\textbf {Theorem 2} Forward  or  mixed stability analysis or any post priori error analysis  is  absolutely invalid stability analyses  for any algorithm in the sense that they are absolutely invalid analyses  for detecting potential instabilities of any algorithm due to excessive element growth. Any perturbation theory can not
 detect potential instabilities of an algorithm either.

\vskip 1em

\noindent Proof:   Forward error analysis or forward stability analysis depends on forward error measurement abs(y' - y) which measures how much a point y' moves away from a reference point y, where y' is the computed solution for exact y=f(x).
It has   a wrong reference point y which could  be moving with y' in the same direction during a computation.   It is therefore invalid for  detecting potential instabilities of any  algorithm due to excessive element  growth in y' relative to  x which is fixed during computation for y.  In contrast, backward stability analysis has the right reference point x.  Similar argument can be applied to mixed stability analysis. That is, the mixed error analysis has a wrong reference point too.  Proof done.
\vskip 1em

\noindent In the proof of Theorem 2, it is pointed out that forward or mixed stability analysis is  absolutely invalid stability analysis in the sense that they are absolutely invalid  for detecting potential instabilities due to excessive element growth because they have absolutely wrong reference points. In contrast,   backward error measurement abs(x - x') has a right reference point x. But having a right reference point itself  is not sufficient.   According to the standard teaching, if a backward error bound is derived for an algorithm f, then the algorithm f is said to be stable("ideal"  case scenario) and backward stability analysis  is absolute,  meaning it does not depends anything.  It is not true in general. When y', the computed output for exact y=f(x)  grows excessively large relative to x  during a computation, its inverse image x' defined by y'=f(x') may or may not grow with it relative to x. It is possible that  it may stay bounded relative to x  despite y' grows excessively large relative to x in which case, a backward stability analysis fails to detect potential instabilities due to excessive element growth. One dimension vector  iteration y=f(x)=1/x is such an example. For  an input x0 approaching to  zero, a computed data y'=1/x0(1+u) is arbitrarily large,  where abs(u) is less than machine precision. f inversely transforms y' to x' defined by y'=f(x') so that x' = x0/(1+u) = x0(1+u+$u^2$+...). Therefore, (x'-x0)/x0  = u +O($u^2$). We see  the computed data y' can be arbitrarily large relative to  x0 but it's inverse image x' defined by y'=f(x') relative to x0 stays bounded. Backward stability analysis fails totally to be able to detect excessive element growth of f when x0 approaches to zero.  See the following  Theorem 3.

\vskip 1em
\noindent\textbf{Theorem 3 }  Backward error analysis is not a stability analysis in general in the sense that it can not detect potential instabilities due to excessive element growth. But for a subclass of algorithms  f satisfying inversely increasing condition,  that is, whenever y' is increasing in a norm excessive large  relative to x, where y' is computed solution for exact y=f(x),  then its inverse image x' defined by y'=f(x') is also increasing in the norm excessively large  relative to x, then a  backward error analysis often known as  backward stability analysis or a priori error analysis   is  able to detect potential instabilities of f due to excessive element growth for this sub  class of algorithms.

\vskip 1em
 \noindent Proof: Self evident by the definition of backward error analysis.

\vskip 1em
\noindent\textbf{Corollary  }  QR and Jacobi's method clearly   belong to this subclass of algorithms because they  inversely  preserve  a  distance.
\vskip 1em

\noindent\textbf{Theorem 4 } No error analysis is able to detect potential instabilities of an algorithm due to excessive element growth in general. Nor is any perturbation result.  Whether or not an algorithm has potential instabilities or not can not be determined by deriving  any relative   error bound with or without any perturbation result.

\vskip 1em
\noindent Proof: combining Theorem 2,3 above. Any perturbation result can not be applied to analyze any algorithm to detect and expose potential instabilities of an algorithm.  Proof done.
\vskip 1em

 \noindent\textbf {Theorem 5}:  No one is sure that she/he can compute any digit correctly in any numerical computing   according to standard teaching of today that backward stability analysis is an  absolute stability analysis for a generic algorithm as in \cite{DM},p.5 or  based on forward or mixed stability analysis with or without a perturbation theory.

\vskip 1em
 \noindent  Proof:   Any error analysis is not able to detect potential instabilities due to excessive element growth in general according to Theorem 2, 3 above. Any matrix perturbation can not determine if an algorithm has potential instabilities or not either.  It does not have any information about any algorithm. Therefore, no body is sure if any algorithm has potential instabilities or not according to the standard teaching of today.  And a  testing in \cite{YY}  shows that no correct digit can be expected in a computed solution  in general if instabilities due to excessive element growth occur.

 \vskip 1em
 \noindent In conclusion,  backward stability analysis is not an absolute stability analysis. It depends on algorithms.  It has been taught as absolute stability analysis, just as time is taught as absolute time in classical physics. We have proved that it is impossible to prove that an algorithm is stable according to the standard teaching of today,just as it is impossible to prove mathematical equations(e.g. Maxwell) is a law of physics according to standard teaching in classical physics.  None of backward stable or forward stable  or mixed stable algorithm  as in \cite{DV},\cite{DK} and \cite{FP} is truly  stable in general. It may have potential instabilities due to excessive element growth.

\end{document}